\documentclass[11pt,a4paper]{amsart}
\usepackage[colorlinks=true, pdfstartview=FitV, linkcolor=blue, 
citecolor=blue]{hyperref}

\usepackage{amsmath}
\usepackage{amsthm}
\usepackage{amssymb,amscd}
\usepackage{microtype,xcolor,enumerate,comment}
\usepackage{a4wide}
\usepackage {diagbox} 
\usepackage{float}

\theoremstyle{plain} 
\newtheorem{theorem}{Theorem}
\newtheorem{corollary}{Corollary}
\newtheorem{lemma}{Lemma}

\theoremstyle{definition} 

\theoremstyle{remark} 

\begin{document}

\title[On the behavior of the number of $k$-full integers]
{On the behavior of the number of $k$-full integers between successive $k$-th powers}

\author{Shusei Narumi}
\author{Yohei Tachiya}
\address{Graduate School of Science and Technology,\newline
\hspace*{\parindent}Hirosaki University, Hirosaki 036-8561, Japan}
\email{h25ms125@hirosaki-u.ac.jp}
\email{tachiya@hirosaki-u.ac.jp}

\makeatletter
\@namedef{subjclassname@2020}{%
  \textup{2020} Mathematics Subject Classification}
\makeatother


\date{\today}

\begin{abstract} 
Let $k\ge2$ be an integer. 
Inspired by an observation of Erd\H{o}s (1976), we investigate the behavior of $h_k(n)$, the number of $k$-full integers in the interval $(n^k, (n+1)^k)$ for 
positive integers $n$. 
Let $\delta_k$ denote the supremum of the exponents $\delta$ for which the inequality $h_k(n) > (\log n)^\delta$ holds for infinitely many positive integers $n$. 
We improve the previously known lower bound $\delta_k \ge 1/2$, implied by Munsch, Shparlinski, and Yau (2020), to $\delta_k \ge 1$. 
Furthermore, for any given integer $m \ge 1$, we derive an effective upper bound for the smallest positive integer $n$ satisfying $h_k(n)\ge m$. \\
\begin{flushright}
-- {\it Dedicated to the memory of our friend Florian Luca.}
\end{flushright}
\end{abstract}

\maketitle

\section{Introduction and Main results}\label{sec1}

Throughout this paper, let 
$k\geq2$ be an integer. 
A positive integer $n$ is called a $k$-full integer 
if $p^k$ divides $n$ for every prime factor $p$ of $n$. 
When $k=2$, such integers are referred to as square-full or powerful integers. 
On the other hand, an integer $n$ is called square-free if it is not divisible by the square of any prime. 
It is well known that the set $\mathcal{F}_k$ of all $k$-full integers can be expressed as the disjoint union  
\begin{equation}\label{eq:10289}
\mathcal{F}_k = \bigsqcup_{\lambda \in \{1\}\cup \Lambda_k} \{x^k \lambda^k \mid x \in \mathbb{Z}_{\geq 1}\},
\end{equation}
where the set $\Lambda_k$ is defined by 
\begin{equation}\label{eq:lambda_def}
\Lambda_k:=\left\{
(y_1^{k+1}\cdots y_{k-1}^{2k-1})^{\frac{1}{k}}\,\middle|\, 
\substack{y_1,\ldots,y_{k-1} \in \mathbb{Z}_{\geq 1},\,\,\,
y_1\cdots y_{k-1}\ge2,\\y_1\cdots y_{k-1} \text{is square-free}}
\right\}.
\end{equation}

In what follows, let $h_k(n)$ denote 
the number of $k$-full integers in the interval $(n^k, (n+1)^k)$ for 
positive integers $n$, 
namely,
\begin{equation}\label{eq:1027111}
h_k(n):=\#\{m\in\mathcal{F}_k\mid n^k<m<(n+1)^k\},\quad n\ge1.
\end{equation}
In \cite[p.~35]{E75}, Erd\H{o}s mentioned without proof that 
\begin{equation}\label{eq:E75_sup}
\limsup_{n \to \infty} h_2(n)=\infty,
\end{equation}
and that, for every integer $\ell \geq 0$, there exists a constant $d_\ell\ge 0$ such that  
\begin{equation}\label{eq:E75_dens}
\#\{ n\in\mathbb{Z}_{\ge1} \mid h_2(n) = \ell, \,1 \leq n \leq x\}
=d_\ell x + o(x) \quad (x \to \infty),
\end{equation}
with $\sum_{\ell=0}^\infty d_{\ell}=1$. 
Moreover, he suspected that 
\begin{equation}\label{eq:E75_conj}
\sup \{\delta \geq 0\mid h_2(n)>(\log n)^\delta \ \text{for infinitely many integers $n$}\} <\infty.
\end{equation}
Subsequently, Shiu~\cite{S80} rigorously proved \eqref{eq:E75_dens}, with $d_\ell>0$ for every integer $\ell\ge0$, 
which in particular established \eqref{eq:E75_sup}. 
He also provided the first few explicit values of $d_\ell$, such as $d_0=0.2759\dots$, $d_1=0.3955\dots$, and $d_2=0.2312\dots$. 
Xiong and Zaharescu~\cite{XZ11} generalized Shiu's result to $k$-full integers by showing that, 
for every integer $\ell \geq 0$, 
there exists an explicit positive constant $d_\ell^{(k)}$ such that
\[
\#\{ n\in\mathbb{Z}_{\ge1}\mid h_k(n) = \ell,\,1 \leq n \leq x  \}=
d_\ell^{(k)} x + O(x/(\log \log x)^\frac{1}{2k}) \quad (x \to \infty),
\]
and moreover, they obtained the generating function of the sequence $(d_{\ell}^{(k)})_{\ell\ge0}$:
\begin{equation}\label{eq:XZ11_genf}
\sum_{\ell=0}^\infty d_{\ell}^{(k)}z^\ell = \prod_{\lambda \in \Lambda_k} \left(1+\frac{z-1}{\lambda}\right).
\end{equation}
Substituting $z=1$ into \eqref{eq:XZ11_genf} yields $\sum_{\ell=0}^{\infty}d_\ell^{(k)}=1$, as Erd\H{o}s observed for $k=2$. 
Recently, the authors~\cite{NT26} investigated the distribution of $k$-full integers in the wider intervals $(n^k,(n+2)^k)$ and established that, 
for any integers $\ell,m\ge0$, the set 
\[
\{ n\in\mathbb{Z}_{\ge1}\mid h_k(n) = \ell,\,h_k(n+1) = m\}
\] 
has positive asymptotic density.    

On the other hand, 
the assertion \eqref{eq:E75_conj} remains open, 
although some partial results have been obtained. 
Define 
\[
\delta_k := \sup 
\big\{\delta \geq 0 \mid h_{k}(n) > (\log n)^\delta \ \text{for infinitely many integers} \ n \big\}
\in 
\mathbb{R}_{\geq 0} \cup \{\infty\}.
\]
The assertion \eqref{eq:E75_conj} can then be restated as $\delta_2 < \infty$. 
De Koninck and Luca~\cite{DL04} proved that, for $k=2$, there exist infinitely many positive integers $n$ such that
$h_2(n) > \frac{9}{20} (\frac{\log n}{\log \log n})^{\frac{1}{3}}$.
This result was subsequently generalized to $k$-full integers and improved by De Koninck, Luca, and Shparlinski~\cite{DLS05}, who showed that
$h_k(n) \geq ( (\frac{3}{8} + o(1)) \frac{\log n}{\log \log n} )^{\frac{1}{3}}$
for infinitely many positive integers $n$.
Moreover, Munsch, Shparlinski, and Yau~\cite{MSY20} further improved this lower bound to
\[
h_{k}(n) \geq \left( \Big(\frac{k}{2(k+1)} + o(1)\Big)\frac{\log n}{\log \log n} \right)^{\frac{1}{2}},
\]
thereby establishing $\delta_k \geq 1/2$ for every integer $k \ge 2$.\\

The aim of this paper is to give new lower bounds for $h_k(n)$ and $\delta_k$. 
Our main results are the following. 

\begin{theorem}\label{thm1}
Let $k\ge2$ be an integer and let $\alpha$ be a constant satisfying  
\begin{equation}\label{eq:10299}
0<\alpha\le \frac{k}{2(k+1)}. 
\end{equation}
Then, as $x\to\infty$, we have 
\begin{align*}
& \log \#\left\{ n\in\mathbb{Z}_{\ge1} \,\middle|\, h_{k}(n)>\alpha\frac{\log x}{\log\log x},\,1 \leq n \leq x \right\} \\
& \quad \ge \left(1-2\alpha\Big(1+\frac{1}{k}\Big)\right)\log x + 2\alpha\Big(1+\frac{1}{k}\Big)\frac{\log x}{\log\log x}\log\log\log x + O\left(\frac{\log x}{\log\log x}\right).
\end{align*}
\end{theorem}

In particular, taking $\alpha =\frac{k}{2(k+1)}$ in Theorem~\ref{thm1} yields   

\begin{corollary}\label{cor1}
Let $k\ge2$ be an integer. 
Then there exist infinitely many positive integers $n$ such that 
\[
h_{k}(n) > \frac{k}{2(k+1)}\cdot\frac{\log n}{\log\log n}.
\]
\end{corollary}

Corollary~\ref{cor1} improves the previously known lower bounds for $h_k(n)$ in 
\cite{DL04,DLS05,MSY20}, thereby establishing a new lower bound for $\delta_k$:  

\begin{corollary}\label{cor2}
For every integer $k\ge2$, we have $\delta_k \geq 1$. 
\end{corollary}

For an integer $m\ge1$, we define 
\[
\eta_{k}(m):=\min\big\{n\in\mathbb{Z}_{\ge1}\mid h_k(n)\ge m\big\}.
\]
By \eqref{eq:10289} and \eqref{eq:lambda_def}, there are no $k$-full integers between $1$ and $2^k$, whereas 
$2^{k+1}$ lies between $2^k$ and $3^k$. 
Hence, $h_k(1)=0$ and $h_k(2)\ge1$, so that $\eta_k(1)=2$ for every integer $k\ge2$. 
A similar argument, together with an examination of the first few square-full and cube-full integers (cf. OEIS A001694 and A036966), 
gives 
\[
\eta_2(2)=5,\quad \eta_3(2)=6,\quad \eta_k(2)=2\,\,\, (k\ge 4).
\]
The following theorem establishes explicit upper bounds for $\eta_{k}(m)$: 
\begin{theorem}\label{thm2}
Let $k\ge2$ be an integer. 
For every integer $m\ge1$, we have  
\[
\log \eta_k(m)
\le 
\frac{k+1}{k}\left(2m\log m +2\Big(\log\frac{\pi^2}{3}-1\Big)m+\sqrt{m}\right).
\]
\end{theorem}

The proofs of Theorems~\ref{thm1} and \ref{thm2} are given in Sections~\ref{sec2} and \ref{sec3}, respectively. 
It should be noted that, instead of relying on Roth's theorem as in \cite{DLS05,MSY20}, 
our proofs employ Lemma~\ref{lem:dirichlet}, an effective refinement of Dirichlet's approximation theorem. 

\section{Proof of Theorem~\ref{thm1}}\label{sec2}
For a real number $x$, $\lfloor x\rfloor$ denotes the greatest integer not exceeding $x$. 
Moreover, the fractional part of $x$ and the distance from $x$ to the nearest integer are defined by 
$\{x\}:=x-\lfloor x \rfloor$ and 
$\|x\|:=\min \{\{x\}, 1-\{x\} \}$, respectively. 
To prove Theorems~\ref{thm1} and \ref{thm2}, we need two lemmas. 
Lemma~\ref{lem:dirichlet} is used in the proof of Lemma~\ref{lem:222}. 
Although the proof of Lemma~\ref{lem:dirichlet} parallels that of the classical Dirichlet approximation theorem, we include it for completeness.

\begin{lemma}\label{lem:dirichlet}
Let $m,n\ge1$ be integers, $\xi_1, \ldots, \xi_m$ be real numbers, and $Q_1,\ldots, Q_m$ be positive integers. 
Then there exist at least $n$ integers $q$ satisfying $1\leq q \leq n Q_1\cdots Q_m$ and    
\[
\|\xi_j q\|<\frac{1}{Q_j},\quad  j =1,\ldots, m. 
\]
Moreover, if $\|\xi_1\|\ge 1/Q_1$, then such integers $q$ can be chosen so that no two of them are consecutive. 
\end{lemma}
\begin{proof}
For $\ell = 1, 2,\ldots, 1+nQ_1\cdots Q_m$, we define $m$-tuples $\mathbf{t}_\ell$ in $[0,1)^m$ by
\[
\mathbf{t}_\ell:=(\{\xi_1 \ell\},\ldots,\{\xi_m \ell\}). 
\]
Consider the partition $\mathcal{P}$ of $[0,1)^m$ into $Q_1 \cdots Q_m$ sub-rectangles, defined by
\[
\mathcal{P} := \left\{ \prod_{j=1}^m \left[ \frac{k_j-1}{Q_j}, \frac{k_j}{Q_j} \right) \middle|\, k_j=1,\dots, Q_j \right\}.
\]
By the pigeonhole principle, there exists an $I\in \mathcal{P}$ containing at least $n+1$ $m$-tuples 
among $\mathbf{t}_1,\ldots,\mathbf{t}_{1+nQ_1\cdots Q_m}$. 
Let 
\begin{equation}\label{eq:1109123}
\mathbf{t}_{\ell_1},\ldots, \mathbf{t}_{\ell_{n+1}}\quad (\ell_1 <\dots < \ell_{n+1})
\end{equation}
be $m$-tuples that lie in $I$. 
Then, for each integer $j=1,\ldots, m$, we have  
\[
\left|\{\xi_j \ell_i\}-\{\xi_j \ell_1\}\right|<\frac{1}{Q_j}, \quad i=2,\ldots,n+1,
\]
and so $\|\xi_j q_i\|<1/Q_j$ for the $n$ distinct positive integers $q_i:=\ell_i-\ell_1$ $(i=2,\ldots, n+1)$ with  
$1\leq q_i\leq nQ_1\cdots Q_m$. 

Next, we prove the latter assertion. Assume $\|\xi_1\|\ge 1/Q_1$. 
Suppose, to the contrary, that there exist two consecutive integers $\ell$ and $\ell+1$ in \eqref{eq:1109123} 
such that $\mathbf{t}_{\ell},\mathbf{t}_{\ell+1}\in I$.   
Then 
\[
\left|\{\xi_1 (\ell+1)\}-\{\xi_1\ell \}\right|<\frac{1}{Q_1}.
\]
On the other hand,  
\[
\left| \{\xi_1(\ell+1)\} - \{\xi_1 \ell\} \right| = \begin{cases}
\{\xi_1\} & \text{if}\,\,\,\{\xi_1 \ell\} + \{\xi_1\} < 1, \\
1 - \{\xi_1\} & \text{otherwise}
\end{cases}
\ge \|\xi_1\|\ge \frac{1}{Q_1},
\]
which is a contradiction. 
Hence, no two integers $\ell$ in \eqref{eq:1109123} are consecutive, 
and therefore, no two of the integers $q_i:=\ell_i-\ell_1$ $(i=2,\dots,n+1)$ are 
consecutive either.
Lemma~\ref{lem:dirichlet} is proved.
\end{proof}
Let $\Lambda_k$ be as defined in \eqref{eq:lambda_def}, and write $\Lambda_k=\{\lambda_j\mid 2<\lambda_1<\lambda_2<\cdots\}$. 
Recall from \eqref{eq:1027111} that $h_k(n)$ denotes the number of $k$-full integers in $(n^k,(n+1)^k)$.
\begin{lemma}\label{lem:222}
Let $k\geq 2$ and $m, n\geq1$ be integers. 
Then there exist at least $n$ integers $q$, no two of which are consecutive, such that 
$2\leq q \leq n\textstyle{\prod_{j=1}^m} (1+\lfloor \lambda_j \rfloor)$ and    
\[
h_k(q-1)+h_{k}(q) \geq m.
\]
\end{lemma}
\begin{proof}
Let $\xi_j:=1/\lambda_j$ and $Q_j:=1+\lfloor\lambda_j\rfloor$ $(j=1,\dots,m)$. 
Then we have $\|\xi_1\|=1/\lambda_1\ge 1/Q_1$. 
Hence, applying the latter assertion of Lemma~\ref{lem:dirichlet}, 
we find that there exist at least $n$ integers $q_1,\ldots, q_n$ $(q_1<\cdots<q_n)$ with $q_{i+1}-q_i\ge 2$ $(i=1,\dots,n-1)$ satisfying 
\[
1\leq q_i \leq n\prod_{j=1}^m (1+\lfloor \lambda_j \rfloor), \quad i=1,\ldots, n,
\]
and 
\begin{equation}\label{ineq:lambda_approx}
\left\|\frac{1}{\lambda_j}q_i\right\|<\frac{1}{1+\lfloor \lambda_j \rfloor}, \quad  i=1,\ldots, n,\quad j=1,\ldots m.
\end{equation}
Inequality \eqref{ineq:lambda_approx} implies that, for each $i=1,\ldots, n$ and $j=1,\ldots, m$, 
there exists an integer $\ell_{i,j}$ such that  
\[
\left|\frac{1}{\lambda_j}q_i-\ell_{i,j}\right|<\frac{1}{\lambda_j},
\]
which shows that all $\ell_{i,j}$ are positive and, since $\lambda_j>2$, that $q_i\ge2$ for each $i$. 
Thus, we have 
\begin{equation}\label{08311}
(q_i-1)^k<\ell_{i,j}^k \lambda_j^k<(q_i+1)^k, \quad  i=1,\ldots, n,\quad j=1,\ldots, m,
\end{equation}
where the integers $\ell_{i,j}^k \lambda_j^k$ are $k$-full integers, but not perfect $k$-th powers. 
Moreover, for a fixed integer $i$, the integers $\ell_{i,j}^k \lambda_j^k$ $(j=1,\dots,m)$ are distinct by \eqref{eq:10289}. 
Hence, by \eqref{08311}, we obtain $h_{k}(q_i-1)+h_k(q_i) \ge m$ for each $i=1,\dots,n$ 
corresponding to the $n$ distinct integers $q_1,\ldots, q_n$.  
The proof of Lemma~\ref{lem:222} is completed. 
\end{proof}
Let $s_j$ be the $j$-th square-free integer. 
Then $s_{j+1}^{1+1/k}\in \Lambda_k$ $(j\geq1)$ by definition~\eqref{eq:lambda_def}. 
Hence, by using the well-known inequality $s_j\le 2j$ $(j\ge1)$, we have    
\begin{equation}\label{ineq:lem2_lambda}
1+\lambda_j\le 1+s_{j+1}^{1+1/k}\le s_{j+2}^{1+1/k}\le \big(2(j+2)\big)^{1+1/k}
\end{equation}
for every integer $j\ge1$, and so 
\begin{align}
\sum_{j=1}^{n}\log(1+\lambda_j)& \leq \left(1+\frac{1}{k}\right)\big(n\log2+\log(n+2)!\big)\nonumber\\
&\le \left(1+\frac{1}{k}\right)n\log n+O(n)\quad (n\to\infty).\label{eq:08312}
\end{align}
\begin{proof}[Proof of Theorem~\ref{thm1}]
Let $\alpha$ be a positive constant as in \eqref{eq:10299}. 
Let $x>0$ be sufficiently large, and define the integers $M:=M(x)$ and $N:=N(x)$ by 
\begin{equation}\label{eq:10295}
M:=\left\lfloor \alpha\frac{\log x}{\log\log x} \right\rfloor+1, \qquad 
N:=\left\lfloor 
\frac{x}{\prod_{j=1}^{2M-1} (1+\lfloor \lambda_j\rfloor)}\right\rfloor.
\end{equation}
Note that the integer $N$ is positive; indeed, by \eqref{eq:08312}, we have 
\begin{align*}
\log \prod_{j=1}^{2M-1} (1+\lfloor \lambda_j\rfloor)&\le 
\sum_{j=1}^{2M-1} \log (1+\lambda_j)\le 
\left(1+\frac{1}{k}\right)2M\log(2M)+O(M)\\
&\le 2\alpha\left(1+\frac{1}{k}\right)\left(1-\frac{\log\log\log x}{\log\log x}\right)\log x+O\left(\frac{\log x}{\log\log x}\right)\quad (x\to\infty),
\end{align*}
and so, from \eqref{eq:10295},  
\begin{align}
\log(N+1)&\ge
\log x-\log \prod_{j=1}^{2M-1} (1+\lfloor \lambda_j\rfloor)\nonumber\\
&\ge \left(1-2\alpha\Big(1+\frac{1}{k}\Big)\right)\log x
+2\alpha\left(1+\frac{1}{k}\right)\frac{\log x}{\log\log x}\log\log\log x+
O\left(\frac{\log x}{\log\log x}\right)\label{eq:10297}
\end{align}
as $x\to\infty$. 
Inequality \eqref{eq:10297} with \eqref{eq:10299} yields $N \to \infty$ as $x\to\infty$. 

Thus, applying Lemma~\ref{lem:222} with $m:=2M-1$ and $n:=N$, we find that 
there exist $N$ positive integers $q_1,\ldots,q_N$ 
($q_1<\cdots<q_N$) with $q_{i+1}-q_i\ge2$ $(i=1,\dots,N-1)$ satisfying 
\[
2\leq q_i \le N\prod_{j=1}^{2M-1}(1+\lfloor \lambda_j \rfloor)\le x
\]
and 
\[
h_k(q_i-1)+h_{k}(q_i) \geq 2M-1, \quad i=1,\dots,N;
\]
in particular, at least one of $h_k(q_i-1)\ge M$ and $h_{k}(q_i) \geq M$ holds for each $i$. 
Define 
\[
n_i := \begin{cases} 
q_i & \text{if $h_k(q_i) \geq M$}, \\ 
q_i-1 & \text{otherwise},
\end{cases} 
\quad i=1, \ldots, N,
\]
where the integers $n_1,\ldots,n_N$ are distinct since $q_{i+1}-q_i\ge2$ $(i=1,\dots,N-1)$. 
Hence, for the set   
\[
\mathcal{S}:=\{n\in\mathbb{Z}_{\ge1}\mid h_{k}(n)\ge M,\,1\leq n \le x \},
\]
we have $\#\mathcal{S}\ge N$, and therefore, it follows from $M>\alpha\log x/\log\log x$ and \eqref{eq:10297} that  
\begin{align*}
& \log \#\left\{ n\in\mathbb{Z}_{\ge1} \,\middle|\, h_{k}(n)>\alpha\frac{\log x}{\log\log x},\,1 \leq n \leq x \right\}\\
&\quad \ge \log \#\mathcal{S} \\
& \quad \ge \log N = \log(N+1)-\log\left(1+\frac{1}{N}\right) \\
& \quad \ge \left(1-2\alpha\Big(1+\frac{1}{k}\Big)\right)\log x + 2\alpha\left(1+\frac{1}{k}\right)\frac{\log x}{\log\log x}\log\log\log x + O\left(\frac{\log x}{\log\log x}\right)
\end{align*}
as $x\to\infty$. 
The proof of Theorem~\ref{thm1} is completed. 
\end{proof}
\section{Proof of Theorem~\ref{thm2}}~\label{sec3}
Recall that $s_j$ denotes the $j$-th square-free integer and that $\Lambda_k=\{\lambda_j \mid \lambda_1<\lambda_2<\cdots\}$ denotes the set of 
real numbers defined in \eqref{eq:lambda_def}. 
To prove Theorem~\ref{thm2}, 
we require a sharper upper bound for $\sum_{j=1}^{m}\log(1+\lfloor \lambda_j\rfloor)$ than the one obtained in \eqref{eq:08312}. 
\begin{proof}[Proof of Theorem~\ref{thm2}]
Let $m\ge1$ be an integer. 
Applying Lemma~\ref{lem:222} with $n=1$, we find that 
there exists an integer $q$ with $2\leq q \leq \textstyle{\prod_{j=1}^{2m-1}} (1+\lfloor \lambda_j \rfloor)$ such that at least one of $h_k(q-1)\geq m$ and $h_k(q)\geq m$ holds.
Hence, $\eta_k(m)\le \prod_{j=1}^{2m-1}(1+\lfloor \lambda_{j} \rfloor)$ 
and, by \eqref{ineq:lem2_lambda}, we have 
\begin{equation}\label{eq:1103}
\log \eta_k(m)\le \left(1+\frac{1}{k}\right)\sum_{j=3}^{2m+1}\log s_{j}.
\end{equation}

In what follows, we derive an upper bound for the right-hand side of \eqref{eq:1103}. 
Let $Q(x)$ denote the number of square-free integers not exceeding $x$, and 
define   
\[
R(x):=Q(x)-\frac{6}{\pi^2}x,\quad x\ge1.
\]
Then 
\begin{equation}\label{eq:1112}
R(n)\geq -c_1\sqrt{n}
\end{equation}
holds for every integer $n\ge1$, where $c_1:=0.103229$. 
Indeed, by \cite[Theoreme~3 bis]{CDM07}, we have  
$|R(x)|\le 0.02767\sqrt{x}$ for $x\ge 438653$, while a direct computation in Mathematica gives  
\[
\min_{1\le n\le 438652}\frac{R(n)}{\sqrt{n}}=\frac{R(380)}{\sqrt{380}}=-0.103228745\dots.
\]
Thus, substituting $n=s_j$ for $j\ge3$ into \eqref{eq:1112}, we have  
\[
\frac{6}{\pi^2}s_j-c_1\sqrt{s_j}- j\le0.
\]
Solving the above quadratic inequality for $\sqrt{s_j}$ yields  
\[
s_j\le \frac{\pi^2}{6}j\Big(f(j)+\sqrt{1+f(j)^2}\Big)^2,
\]
where $f(j):=\sqrt{c_2/j}$ $(j\ge1)$, with $c_2:=(c_1\pi)^2/24$. 
Hence, taking logarithms and summing over $j=3,\dots,2m+1$, we obtain  
\begin{equation}\label{eq:11031}
\sum_{j=3}^{2m+1}\log s_j\le (2m-1)\log\frac{\pi^2}{6}+\log (2m+1)!-\log2+
2\sum_{j=3}^{2m+1}\log\Big(f(j)+\sqrt{1+f(j)^2}\Big).
\end{equation}
From the inequality 
\[
\log n!<n\log n-n+\frac{1}{2}\log(2\pi n)+\frac{1}{12n}
\]
for integers $n\ge1$ (cf.~\cite{R55}), we have 
\begin{align}
\log (2m+1)!&<(2m+1)\log(2m+1)-(2m+1)+\frac{\log (2\pi(2m+1))}{2}+\frac{1}{12(2m+1)}\nonumber \\
&=2m\log m-2(1-\log2)m+T(m),\label{eq:110399}
\end{align}
where 
\[
T(m):=2m\log\Big(1+\frac{1}{2m}\Big)+\log(2m+1)-1+\frac{\log(2\pi (2m+1))}{2}+\frac{1}{12(2m+1)}.
\]
Moreover, 
\begin{equation}\label{eq:11034}
\sum_{j=3}^{2m+1}\log \Big(f(j)+\sqrt{1+f(j)^2}\Big)\le 
\sqrt{c_2}\sum_{j=3}^{2m+1}j^{-1/2}\le 
2\sqrt{c_2}(\sqrt{2m+1}-\sqrt{2}),
\end{equation}
where the first and second inequalities follow from 
$\log(x+\sqrt{1+x^2})\le x$ $(x\ge 0)$ and an integral comparison, respectively. 
Hence, combining \eqref{eq:11031}, \eqref{eq:110399}, and \eqref{eq:11034}, we obtain 
\begin{equation}
\sum_{j=3}^{2m+1}\log s_j<
2m\log m +2\Big(\log\frac{\pi^2}{3}-1\Big)m+
4\sqrt{c_2}(\sqrt{2m+1}-\sqrt{2})+T(m)-\log \frac{\pi^2}{3}.\label{eq:1104}
\end{equation}
Numerical computations in Mathematica show that 
\begin{equation}\label{eq:1111}
4\sqrt{c_2}(\sqrt{2m+1}-\sqrt{2})+T(m)-\log \frac{\pi^2}{3}<\sqrt{m}
\end{equation}
holds for every integer $m\ge 167$. 
For $1\le m\le 166$, we also verify that 
\begin{equation}\label{eq:11112}
\sum_{j=3}^{2m+1}\log s_j<
2m\log m+2\Big(\log\frac{\pi^2}{3}-1\Big)m+\sqrt{m}.
\end{equation}
Combining \eqref{eq:1103}, \eqref{eq:1104}, and \eqref{eq:1111} for $m\ge 167$, and \eqref{eq:1103} and \eqref{eq:11112} for $1\le m\le 166$, we obtain Theorem~\ref{thm2}.
\end{proof}

\section*{Acknowledgments}
This work was supported by JSPS KAKENHI Grant Numbers JP25K06911.

\end{document}